\documentclass[10pt,oneside,a4paper,reqno]{amsart}%{article}
\usepackage[english]{babel}
\usepackage{amsmath,amsfonts,amscd,amssymb,latexsym}
\usepackage{hyperref}

\usepackage{hyperref}
\newtheorem{Th}%[Th]  
{Theorem}
\newtheorem{Lem}[Th]{Lemma}
\newtheorem{cor}[Th]{Corollary}
\newtheorem{exa}[Th]{Example}

\begin{document}

\sloppy

\begin{center}

{\bf
%Free superalgebras of Jordan brackets\\
Algebras of Jordan brackets and Generalized Poisson algebras
\\
}

\
Ivan Kaygorodov$^{a,b}$\footnote{The author was supported
by RFBR 15-31-21169  and by FAPESP 11/51132-9  and  14/24519-8.}

\

$^{a}$ Universidade Federal do ABC, Santo Andre, Brazil. \\
$^{b}$ Sobolev Institute of Mathematics, Novosibirsk, Russia.

\

 e-mail: kib@ime.usp.br

\end{center}

%\begin{center}\textbf{\large \ Ivan Kaygorodov}\\ \end{center}

%\begin{center} \textbf{\large \ kib@math.nsc.ru }\par \end{center}

\medskip

\begin{center}

{\bf Abstract.}

\end{center}

We construct a basis of free unital generalized Poisson superalgebras
and a basis of free unital superalgebras of Jordan brackets.
Also, we prove the analogue of Farkas' Theorem for PI unital generalized Poisson algebras and 
PI unital algebras of Jordan brackets. 
Relations beetwen generic Poisson superalgebras and superalgebras of Jordan brackets will be studied.

\


\medskip

%\begin{minipage}{400} 
%\tableofcontents 
%\end{minipage}

%\newpage

\section{Introduction.}
Many interesting and important results have been obtained about the structure of polynomial algebras, 
free associative algebras, and free Lie algebras. 
Although free Poisson algebras are closely connected with these algebras, only a few results on their structure are known.
Free Poisson superalgebras were introduced in \cite{Shestakov93}.
For example, the free Poisson algebra of rank three and Poisson brackets were used recently in \cite{SHU1,SHU2} 
to prove that the Nagata automorphism 
of the polynomial algebra $F[x, y, z]$ over a field $F$ of characteristic 0 is wild.
Systematic study of free Poisson algebras was initiated in \cite{MLU_Centralizers},
and 
subsequently continued in \cite{MLU_Derivations}-\cite{Umi_Universal_env}.
Generic Poisson algebras were studied in \cite{KSML}.
Other generalizations of Poisson superalgebras are generalized Poisson superalgebras and superalgebras of Jordan brackets.
As follows from \cite{Kantor90,Kantor92}, every Poisson superalgebra is a superalgebra of Jordan brackets. 
In \cite{KingMcc92,KingMcc95} 
all identities whose define unital superalgebras of Jordan brackets was described,
and in \cite{Kay10} analogue for non-unital case was obtained.
Superalgebras of Jordan brackets play important role in classification of simple finite-dimensional Jordan superalgebras \cite{MarZel01}.
The Kantor construction gives interesting relations between Novikov-Poisson algebras and Jordan superalgebras \cite{Zakharov}.

\

Recall the notion of {\bf Kantor's Double \cite{Kantor92}.}
Let $A=A_0 \oplus A_1$ be an associative supercommutative superalgebra and $\{ \ , \ \}: A \times A \rightarrow A$~---
super-anticommutative multiplication, which we shall call a bracket.
 Using superalgebra $A$ and bracket $\{, \}$ we can construct superalgebra $K(A)$.
Consider the following direct sum  of spaces $K(A)=A \oplus  A x$,
where $Ax$ is an isomorphic copy of the space $A.$
Let $a$ and $b$ are homogeneous elements in $A$.
Then multiplication $*$  on $J(A, \{, \})$
is defined by the relations
$$a* b=ab,\  a* bx=(ab)x, $$
$$ax* b=(-1)^{|b|}(ab)x, ax * bx = (-1)^{|b|}\{a, b\}.$$
We set $K(A)_0=A_0 \oplus A_1 x, K(A)_1 = A_1 \oplus A_0 x,$ so  $K(A)$ is a $\mathbb{Z}_2$-graded algebra.

\

That is, superalgebra $JB$ is called a superalgebra of Jordan brackets, if the Kantor Double is a Jordan superalgebra.
As follows from \cite{KingMcc92,KingMcc95}, 
every superalgebra  $JB$ with associative-supercommutative multiplication $\cdot$ and superanticommutative  multiplication $\{,\}$
satisfies 
\begin{eqnarray}\label{KMtojd1}\{a,bc\}=\{a,b\}c+(-1)^{|a||b|}b\{a,c\}-D(a)bc, \end{eqnarray}
$$\{a,\{b,c\}\}=\{\{a,b\},c\}+(-1)^{|a||b|}\{b,\{a,c\}\}+$$
\begin{eqnarray}\label{KMtojd2}D(a)\{b,c\}+(-1)^{|a||bc|}D(b)\{c,a\}+(-1)^{|c||ab|}D(c)\{a,b\},\end{eqnarray}
where $D(a)=\{a,1\}$,
is a superalgebra of Jordan brackets. 
Note that $D$ is an even derivation of superalgebra $(JB,\cdot)$.
If $D=0$, then superalgebra $(JB,\cdot,\{,\})$ is a Poisson superalgebra \cite{Kantor90}.

V. Kac and N. Cantarini studied linearly-compact simple superalgebras of Jordan brackets \cite{KacKant07},
I. Kaygorodov and V. Zhelyabin studied $\delta$-superderivations of simple superalgebras of Jordan brackets \cite{ZheKa11},
the relation between Jordan brackets and Poisson brackets was studied in \cite{marsheze01} by E. Zelmanov, I. Shestakov and C. Martinez, 
and in \cite{ZheKa11,shestakov09}  special superalgebras of Jordan brackets were considered.
As follows from \cite{KacKant07}, we can consider new multiplication 
$\{a,b\}_D= \{a,b\}+\frac{1}{2}(aD(b)-D(a)b)$ and $D_D=\frac{1}{2}D,$ 
and we obtain
new superalgebra with following identities
\begin{eqnarray}\label{jo1}\{a,bc\}_D=\{a,b\}_Dc+(-1)^{|a||b|}b\{a,c\}_D-D_D(a)bc, \end{eqnarray}
\begin{eqnarray}\label{jo2}\{a,\{b,c\}_D\}_D=\{\{a,b\}_D,c\}_D+(-1)^{|a||b|}\{b,\{a,c\}_D\}_D.\end{eqnarray}
It is a generalized Poisson superalgebra.

All vector spaces are considered over fields with characteristic zero.

%\newpage

\section{A basis of free Generalized Poisson superalgebra.}

The description of a basis of a free algebraic structure is one of the main problems in modern algebra.
Researchers 
constructed bases for 
free Lie algebras \cite{shirshov62,chibr06},
free Lie superalgebras \cite{chibr_lie},
partially commutative Lie algebras \cite{poroshenko11} and so on.

Let us now consider free unital generalized Poisson superalgebras.
Let $D_X$ be a free non-associative super-dialgebra over field $k$ with 
set of even and odd generators $X=X_0 \cup X_1=\{1,x_1, \ldots, x_n\}$,
and two multiplications $\cdot$ and $\{,\}.$

We consider ideal $I$ of super-dialgebra $D_X$, generated by 
$$\big\{ \  \{a,b\}+(-1)^{|a||b|} \{b,a\}, \ a\cdot b-(-1)^{|a||b|}b\cdot a, \ 1 \cdot a-a,$$ 
$$\{a,b\}\cdot c+(-1)^{|a||b|} b\cdot \{a,c\}-D(a)\cdot b\cdot c-\{a,b\cdot c\}, $$
$$\{\{a,b\},c\} +(-1)^{|a||b|}\{b,\{a,c\}\}-\{a,\{b,c\}\} | \ a,b,c \in D_{X_0} \cup  D_{X_1}   \   \big\},$$
where $D(a)=\{a,1\}.$
We say, that super-dialgebra $D_X/I$ is a free unital generalized Poisson superalgebra.

Let $GenP$ be a free unital generalized Poisson superalgebra with a set of generators $X$.
We assume, that $X$ is a linearly ordered set and $1 <x_1 < \ldots < x_n.$ 
All elements from $X$ are good words.
Let we defined all good words with length less than $n$.
The word $w$ with length $n$ in alphabet $X$ is a good word, if:

1) $w=\{u,v\}$, where $u,v$  are good words and $u>v$;

2) if $w=\{\{u_1,u_2\},v\},$ then $u_2\leq v.$

We define the order on good words with length less or equal $n$  such 
that it be consistent with order to correct for existing words, and satisfy the conditions
if word $w$ is a good and $w=\{u,v\},$ then $w>u,v.$

Define the set $M$ as
$$M=\big\{u,\{v,v\} | u, v \mbox{ are good words in alphabet  }X, v \in GenP_1\big\}.$$

We define a linear order on the set $M$.
Let $u>v$, if $deg(u)>deg(v)$;
if $deg(u)=deg(v)$, 
then for $u=\{u_1,u_2\}, v=\{v_1,v_2\}$ (where $u_1>u_2, v_1>v_2$) 
we say $u>v$ if $u_1>v_1$, 
and if $u_1=v_1$ we say $u>v$ if $u_2>v_2$.
The elements of set $M$ we denote as $e_i$, where  $e_i<e_j$ if $i<j$.

Let $L_X$ be a free non-associative superalgebra over field $k$ with set of generators $X$ and multiplication $\{,\}_L$.
Consider the ideal $J$ of superalgebra $L_X$, 
generated by 
$$\big\{ \  \{a,b\}+(-1)^{|a||b|} \{b,a\},$$ 
$$\{\{a,b\},c\} +(-1)^{|a||b|}\{b,\{a,c\}\}-\{a,\{b,c\}\} | \ a,b,c \in L_{X_0} \cup  L_{X_1}   \   \big\}.$$
Superalgebra $L_X/J$ is a free Lie superalgebra \cite{shtern}.
As proved in \cite{shtern}, the set $M$ is a basis of the free Lie superalgebra $L_X/J$.

Let $K$ be an ideal in $D_X$ generated by elements with type $\{ x\cdot y | x,y \neq 1\}$. 
It is easy to see, that $(D_X/I)/K  \approxeq L_X/J$.

Let 
$$U=\{ e_{i_1}^{k_1} \ldots e_{i_n}^{k_n} | e_{i_1} \neq 1, e_{i_m} \in M,  i_1 < \ldots < i_n \} \cup \{1\}.$$
The set $U$ is closed under multiplication $\cdot$. 
The multiplication $\cdot$ is a supercommutative,
and we can consider elements of $M$ as 
$e_{i_1}^{k_1} \ldots e_{i_n}^{k_n},$ where $k_i \neq 0,$ and if $e_{i_j} \in GenP_1$, then $k_j=1.$
We say that element $e_{i_1}^{k_1} \ldots e_{i_n}^{k_n}$ has degree $deg_e=\sum k_i$.
For element $e_i \in M$ we define $deg_L$ as the length of word $e_i$ with generators $X$.
The degree of element $f$ is the number of generator elements in word $f$.

The main goal of this chapter is a description of the basis of free unital generalized Poisson superalgebras. 
We will prove that set $U$ is a basis of free unital generalized Poisson superalgebras.

\begin{Lem}\label{vectorspace}
The set $U$ is a basis of free superalgebra $GenP$ as a vector space.
\end{Lem}

{\bf Proof.}
We consider word $w$ with degree $2$. 
It is easy to see that $w$ has type $x_ix_j$, or $\{x_i,x_j\}$, or $\{1,x_i\}$, and $w \in U$.
We prove that every word with degree $n$ is a linear combination of elements from $U$.
That is every element of superalgebra $GenP$ is a linear combination of elements from $U$.
We will use induction.
Let every word with degree less than $n$ be a linear combination of elements from $U.$
Then $w=w_1w_2$ or $w=\{w_1,w_2\}$, where $w_1,w_2$ are words with degree less than  $n$, 
that is, $w_1$ and $w_2$   are linear combinations of elements from  $U$.

First case:
$w_1=\sum_i \alpha_{1i} m_{1i}$ and $w_2=\sum_i \alpha_{2i} m_{2i},$ where $m_{1i},m_{2i} \in U, \alpha_{1i}, \alpha_{2i} \in k.$
That is, $w=\sum_{i,j} \alpha_{1i}\alpha_{2j}  m_{1i}m_{2j}$
and using closed $U$ under multiplication $\cdot$,
we can say that $w$ is a $\sum_i \alpha_i m_i, m_i \in U, \alpha_i \in k.$

Second case:
$w_1=\sum_i \alpha_{1i} m_{1i}$ and $w_2=\sum_i \alpha_{2i} m_{2i},$ where $m_{1i},m_{2i} \in U,\alpha_{1i}, \alpha_{2i} \in k.$
That is,
$w=\sum_{i,j} \alpha_i \alpha_j \{m_{i},m_{j}\}$.
We consider arbitrary summand $\{m_i,m_j\}$. 
If we can consider $m_i$ or $m_j$ as $n_1n_2,$ where $n_i \in U, n_i \neq 1,$
then, we can say that $m_j$ has type $n_1n_2$.
Using the relation
$$\{m_i,n_1n_2\}= \{m_i,n_1\}n_2+(-1)^{|m_i||n_1|}n_1\{m_i,n_2\}-D(m_i)n_1n_2,$$
and note that degrees of elements $\{m_i,n_1\},\{m_i,n_2\},D(m_i)$ less than $n$ and $n_1,n_2 \in U,$
we can say that elements $\{m_i,n_1\},\{m_i,n_2\},D(m_i)$  are linear combinations of elements from $U$,
and $w$ is a linear combination of elements from $U$.

If $m_i$ and $m_j$ is not of type $n_1n_2,$ where $n_1,n_2 \in U$ and $n_1,n_2 \neq 1$, 
then $m_i,m_j \in M.$
That is, if $m_i=m_j \in GenP_1,$ then $\{m_i,m_j\} \in M.$
In the other case, 
$m_i \neq m_j$ and we can say that $m_i>m_j$ and if $m_i=\{m_{1i},m_{2i}\}, m_{1i}\geq m_{2i}$ 
then
for $m_{2i}\leq m_j$ it follows that $\{m_i,m_j\} \in M,$
for $m_{2i}>m_j$ it follows that
$$\{m_i,m_j\} = \{\{m_{1i},m_j\},m_{2i}\}+(-1)^{|m_{1i}||m_j|}\{m_{1i},\{m_{2i},m_j\}\},$$
where $\{\{m_{1i},m_j\},m_{2i}\} \in M.$
If $\{m_{1i},\{m_{2i},m_j\}\} \in M$, then we obtained that $\{m_i,m_j\}$ is a linear combination of elements from $U$,
or we can consider $m_{1i}$ as $\{,\}$-multiplication of elements  $m_{11i}$ and $m_{21i}$.
After that, we can continue the process, using identity (\ref{jo2}).
Note that, 
degree of $m_{i}$ is greater than  degree of $m_{1i}$,
degree of $m_{1i}$ is greater than degree of $m_{11i}$, and so on.
This process will be finite.
That is, the lemma is proved.

\medskip

For element $b=\Pi_{k=1}^n e_k^{t_k}$ 
we say $\frac{b}{e_m}=e_1^{t_1} \cdots e_{m-1}^{t_{m-1}}e_m^{t_m-1} e_{m+1}^{t_{m+1}}\cdots e_n^{t_n}.$

It is easy to see that

$$\{a,e_k^{t_k}\}=
\{a,e_k^{t_k-1}\}e_k +(-1)^{|a||e_k^{t_k-1}|}e_k^{t_k-1} \{a,e_k\} -D(a)e_k^{t_k}=$$
$$\{a,e_k^{t_k-2}\}e_k^2+(-1)^{|a||e_k^{t_k-2}|}e_k^{t_k-2} \{a,e_k\}e_k+(-1)^{|a||e_k^{t_k-1}|}e_x^{t_k-1}\{a,e_k\} -2D(a)e_k^{t_k}=$$
$$\ldots$$
$$=t_k \{a,e_k\} e_k^{t_k-1}-(t_k-1)D(a)e_k^{t_k}.$$

Using this notations, we can conclude
$$\{a,b\}=\{a, \Pi_{k=1}^n e_k^{t_k} \}=
\left\{a, \frac{b}{e_n^{t_n}} \right\}e_n^{t_n}+(-1)^{|a|\left|\frac{b}{e_n^{t_n}}\right|} \frac{b}{e_n^{t_n}} \left\{a, e_n^{t_n} \right\} -D(a)b=$$
$$\left\{a, \frac{b}{e_{n-1}^{t_{n-1}}e_n^{t_n}} \right\}e_{n-1}^{t_{n-1}}e_n^{t_n}+
(-1)^{|a|\left|\frac{b}{e_{n-1}^{t_{n-1}}e_n^{t_n}}\right|} \frac{b}{e_{n-1}^{t_{n-1}}e_n^{t_n}} \left\{a, e_{n-1}^{t_{n-1}} \right\} e_n^{t_n} +$$ 
$$(-1)^{|a|\left|\frac{b}{e_n^{t_n}}\right|} \frac{b}{e_n^{t_n}} \left\{a, e_n^{t_n} \right\} -2D(a)b=$$
$$\ldots$$
$$=\sum_{k=1}^n
(-1)^{|e_1^{t_1}\cdots e_{k-1}^{t_{k-1}}||e_k^{t_k}|}
 t_k\{a, e_k\} \frac{b}{e_k} - (\sum_{k=1}^n t_k-1)D(a)b.$$

\

We define linear mapping $^*:U \rightarrow U$ as:

1) if $w \in X$, then $w^*=w$;

2) if $w \in M$, then $w^*=w$;

3) if $a=\Pi_{k=1}^n e_k^{a_k}, b=\Pi_{k=1}^n e_k^{b_k}$, then
$$(ab)^*= 
(-1)^{\sum_{i=1}^n |e_i^{a_i}| |e_{i+1}^{b_{i+1}}\cdots e_n^{b_n}| } \Pi_{k=1}^n e_k^{a_k+b_k}.$$

4) if we defined mapping $^*$ for all words with length less than $n$, 
then word $w=\{w_1,w_2\}$, where $w_1,w_2 \in M$, 
we can consider as element of free Lie superalgebra $(M, \{,\})$ 
and can define $w^*$ as a linear combination of basis elements of superalgebra $(M, \{,\})$, 
that is $w$ be a linear combination of elements of $M$.

5) if we defined mapping $^*$ for all words with length less than $m$, 
then for word $w=\{w_1,w_2\}$ with length $m$,
where $w_1$ or $w_2$ is an element of $U \setminus M$, 
we can say that $w_2 \in U \setminus M$ and $w_2=\Pi_{k=1}^n e_k^{t_k}$,
then

$(\{w_1,w_2\})^*=$
$$\sum_{k=1}^n 
(-1)^{|e_1^{t_1}\cdots e_{k-1}^{t_{k-1}}||e_k^{t_k}|}
 t_k\left(\{w_1, e_k\}^* \frac{w_2}{e_k}\right)^* - \left(\sum_{k=1}^n t_k-1 \right)(D(w_1)^*w_2)^*,$$
where in the right part of expression elements with type $\{x,y\}^*$ was defined above.

\

%\newpage

\begin{Th}\label{basis}
The set $U$ is a basis of free superalgebra $GenP.$ 
\end{Th}

{\bf Proof.}
We consider set $U$ as a basis of a vector space and define two new multiplications  $*$ and $\{,\}^*$ 
by:
$u*v=(uv)^*$ and $\{u,v\}^*=(\{u,v\})^*.$
It is a superalgebra $(U,*,\{,\}^*).$

We want to prove that superalgebra $(U, *, \{,\}^*)$ is a generalized Poisson superalgebra.
After that we can say that $U$ is a basis of free generalized Poisson superalgebra $GenP.$

It is easy to see
that multiplication $*$ is associative and supercommutative, 
and multiplication $\{,\}^*$ is superanticommutative.
We want to prove that superalgebra $(U,*, \{,\}^*)$ satisfies
\begin{eqnarray}\label{jo1*}\{a,b*c\}^*=\{a,b\}^**c+(-1)^{|a||b|}b*\{a,c\}^*-D(a)^**b*c.\end{eqnarray}
We say that $a,b,c \in U$ and will prove by induction on $deg_e(a)+deg_e(b)+deg_e(c)$.

Note that
$$ \{1,  e_i *e_j \} ^* = \{1, e_i\}^**e_j+e_i* \{1,e_j\}^*,$$
$$\{e_k, e_i *e_j \}^*= \{e_k,e_i\}^**e_j +(-1)^{|e_k||e_i|} e_i* \{e_k, e_j\}^* -D(e_k)^** e_i*e_j$$
by the definition of mapping $^*$. 
That is, 
if $deg_e(a)+deg_e(b)+deg_e(c)<4,$ then identity (\ref{jo1*}) is true.

Let $deg_e(a)+deg_e(b)+deg_e(c)>3.$ 
We can say that $b=\Pi_{k=1}^n e_k^{b_k}, c=\Pi_{k=1}^n e_k^{c_k}$ and 
define new notations:
$$\gamma(b,c)= \sum_{i=1}^n |e_i^{b_i}| |e_{i+1}^{c_{i+1}} \cdots e_{n}^{c_n}| \mbox{ and }
\mu(b,c,k)=|e_1^{b_1+c_1}\cdots e_{k-1}^{b_{k-1}+c_{k-1}}||e_k^{b_k+c_k}|.$$
That is,

\begin{eqnarray}\label{jordok0}
D(a)*b*c= (-1)^{\gamma(b,c)} D(a)*\Pi_{k=1}^n e_k^{b_k+c_k}.
\end{eqnarray}

We can see that

\begin{eqnarray}\label{jordok1}
(-1)^{\gamma(b,c)}\{a,b*c\}^* =
\end{eqnarray}

\begin{eqnarray*}
\sum_{k=1}^n 
(-1)^{\mu(b,c,k)}
 (b_k+c_k) \{a, e_k\}^** \frac{ \Pi_{m=1}^n e_m^{b_m+c_m}}{e_k} - \\
 \left(\sum_{k=1}^n (b_k+c_k)-1\right)D(a)^**  \Pi_{m=1}^n e_m^{b_m+c_m}.
\end{eqnarray*}

It is easy to see that

$$\{a,b\}^**c =
\sum_{k=1}^n  
(-1)^{|e_1^{b_1}\cdots e_{k-1}^{b_{k-1}}||e_k^{b_k}|}
 b_k\{a, e_k\}^**\frac{b}{e_k} *c - \left(\sum_{k=1}^n b_k-1\right)D(a)^* *  b * c=$$
$$\sum_{k=1}^n  
(-1)^{|e_1^{b_1}\cdots e_{k-1}^{b_{k-1}}||e_k^{b_k}|
+ \sum_{i=1}^n |e_i^{b_i}| |e_{i+1}^{c_{i+1}}\cdots e_n^{c_n}|
-|e_k||e_1^{c_1} \cdots e_{k-1}^{c_{k-1} }| }
 b_k\{a, e_k\}^**\frac{\Pi_{m=1}^n e_m^{b_m+c_m}}{e_k}$$
$$ - (-1)^{\sum_{i=1}^n |e_i^{b_i}| |e_{i+1}^{c_{i+1}}\cdots e_n^{c_n}|} \left(\sum_{k=1}^n b_k-1\right)D(a)^* *  \Pi_{m=1}^n e_m^{b_m+c_m}.$$

That is,

\begin{eqnarray}\label{jordok2}
\{a,b\}^**c = 
\end{eqnarray}
$$(-1)^{\gamma(b,c)} 
\left(\sum_{k=1}^n 
(-1)^{\mu(b,c,k)}
 b_k \{a, e_k\}^** \frac{ \Pi_{m=1}^n e_m^{b_m+c_m}}{e_k} - \left(\sum_{k=1}^n b_k-1\right)D(a)^**  \Pi_{m=1}^n e_m^{b_m+c_m}\right). $$

Note that,

$$(-1)^{|c||b|}\{a,c\}^**b =
(-1)^{|c||b|}\sum_{k=1}^n 
(-1)^{|e_1^{c_1}\cdots e_{k-1}^{c_{k-1}}||e_k^{c_k}|}
c_k\{a, e_k\}^** \frac{c}{e_k}*b - \left(\sum_{k=1}^n c_k-1\right)D(a)^* *b*c=$$
$$\sum_{k=1}^n  
(-1)^{ |e_k^{c_k}||b|+|e_k||e_1^{c_1} \cdots e_{k-1}^{c_{k-1}}| +  \sum_{i=1}^n |e_i^{b_i}| |e_{i+1}^{c_{i+1}}\cdots e_n^{c_n}|  - |e_k||e_{k+1}^{b_{k+1}} \cdots e_n^{b_n}|}  b_k\{a, e_k\}^**\frac{\Pi_{m=1}^n e_m^{b_m+c_m}}{e_k}$$
$$ - (-1)^{\sum_{i=1}^n |e_i^{b_i}| |e_{i+1}^{c_{i+1}}\cdots e_n^{c_n}|} \left(\sum_{k=1}^n c_k-1\right)D(a)^* *  \Pi_{m=1}^n e_m^{b_m+c_m}.$$

That is,
\begin{eqnarray}\label{jordok3}
(-1)^{|c||b|}\{a,c\}^**b =
\end{eqnarray}
$$(-1)^{\gamma(b,c)} 
\left(\sum_{k=1}^n 
(-1)^{\mu(b,c,k)}
 c_k \{a, e_k\}^** \frac{ \Pi_{k=1}^n e_m^{b_m+c_m}}{e_k} - \left(\sum_{k=1}^n c_k-1\right)D(a)^**  \Pi_{m=1}^n e_m^{b_m+c_m}\right). $$

In the end, 
summing (\ref{jordok0}) and (\ref{jordok1}) and subtracting (\ref{jordok2}) and (\ref{jordok3}),
we can obtain that superalgebra $(U,*,\{,\}^*)$ satisfies (\ref{jo1*}).

We want to prove that the superalgebra satisfies the identity

\begin{eqnarray}\label{jo2*}\{a,\{b,c\}^*\}^*=\{\{a,b\}^*,c\}^*+(-1)^{|a||b|}\{b,\{a,c\}^*\}^*.\end{eqnarray}

If $a,b,c \in M$, 
then identity (\ref{jo2*}) is true, 
because $(M, \{,\}^*)$ is a free Lie superalgebra.

Let $a=a_1*a_2,$ where $a_1,a_2 \neq 1$ and $a_1,a_2 \in U$, 
then using (\ref{jo1*}), we can obtain 

$\{\{a_1*a_2,b\}^*,c\}^* +(-1)^{|a||b|} \{b,\{a_1*a_2,c\}^*\}^*=$

$$-(-1)^{|b||a|}\{\{b,a_1\}^**a_2+(-1)^{|b||a_1|}a_1 *\{b,a_2\}^* -D(b)^*a,c\}^*-$$
$$(-1)^{|a|(|b|+|c|)} \{b, \{c,a_1\}^**a_2 +(-1)^{|c||a_1|}a_1*\{c,a_2\}^*-D(c)^**a_1*a_2\}^*=$$

$$(-1)^{|b||a|+|c|(|a|+|b|)}(\{c, \{b,a_1\}^*\}^**a_2+(-1)^{(|b|+|a_1|)|c|}\{b,a_1\}^**\{c,a_2\}^*-D(c)^**\{b,a_1\}^**a_2)+$$
$$(-1)^{|b||a|+|c|(|b|+|a|)+|b||a_1|}(\{c,a_1\}^**\{b,a_2\}^*+(-1)^{|c||a_1|}a_1*\{c,\{b,a_2\}^*\}^*-D(c)^**a_1*\{b,a_2\}^*)+$$
$$(-1)^{|b||a|+|c|(|a|+|b|)}(\{c,D(b)^*\}^**a+(-1)^{|c||b|}D(b)^**\{c,a\}^*-D(c)^**D(b)^**a) - $$
$$(-1)^{|a|(|b|+|c|)}(\{b,\{c,a_1\}^*\}^**a_2 +(-1)^{|b|(|c|+|a_1|)}\{c,a_1\}^**\{b,a_2\}^*-D(b)^**\{c,a_1\}^**a_2)-$$
$$(-1)^{|a|(|b|+|c|)+|c||a|}(\{b,a_1\}^**\{c,a_2\}^*+(-1)^{|b||a_1|}a_1*\{b,\{c,a_2\}^*\}^*-D(b)^**a_1*\{c,a_2\}^*)-$$
$$(-1)^{|a|(|b|+|c|)} (\{b,D(c)^*\}^**a+(-1)^{|b||c|}D(c)^**\{b,a\}^*-D(b)^**D(c)^**a)=$$

$$-(-1)^{|a|(|b|+|c|)}(\{\{b,a_1\}^*,c\}^*+(-1)^{|b||a_1|}\{b,\{c,a_1\}^*\}^*)*a_2$$ 
$$+(-1)^{(|b|+|c|)|a_1|}a_1*(\{\{b,a_2\}^*,c\}^*+(-1)^{|b||a_2|}\{b,\{c,a_2\}^*\}^*)-\{D(b)^*,c\}^**a-\{b,D(c)^*\}^**a=$$

$$-(-1)^{|a||bc|}(\{\{b,c\}^*,a_1\}^** a_2+(-1)^{|bc||a_1|}a_1*\{\{b,c\}^*,a_2\}^*-D(\{b,c\}^*)^**a=\{a_1*a_2, \{b,c\}^*\}^*.$$

If $a \in M$, $b \in U \setminus M$ or $c \in U \setminus M$, 
then designating $a,b,c$,
we can use the given calculations. 
That is, the superalgebra satisfies (\ref{jo2*}) and it follows that the theorem is proved.

\medskip

Let $GenP(n)$ be a space of multi-homogeneous components of free generalized Poisson algebra 
$GenP$ with $n+1$ generators $1,x_1, \ldots, x_n$.
Then 

\begin{Th}\label{dimen}
 $dim(GenP(n))=n\cdot n!.$
\end{Th}

{\bf Proof.}
Note that, 
we consider the vector space of free algebra $GenP[1,x_1, \ldots, x_n]$ as a subspace in the free Poisson algebra
$P[1,x_1, \ldots, x_n]$ with $n+1$ generators. 
It is well-known that the dimension of the space of multi-homogeneous components of the free Poisson algebra with $n+1$ generators is a $(n+1)!$.
Using the description of basis of free generalized Poisson algebra (Theorem \ref{basis}), 
we can define the structure of the basis of the space of multi-homogeneous components. 
That is, 
the basis of the space of multi-homogeneous components of the free generalized Poisson algebra with generators $1,x_1, \ldots, x_n$ 
and the basis of the space of multi-homogeneous components of the free Poisson algebra with generators $1,x_1, \ldots, x_n$,
has distinction only for elements $1 \cdot y$, where $y$ is an arbitrary element from the basis of the 
free Poisson algebra with generators $x_1, \ldots, x_n$.
That is, $dim(GenP(n))= (n+1)!-n!=n\cdot n!.$
The Theorem is proved.

%THE BASIS OF SUPERALGEBRA OF JORDAN BRACKETS

\section{A basis of free superalgebra of Jordan brackets.}

Let us now consider free superalgebras of Jordan brackets.
Let $D_X$ be a free non-associative super-dialgebra over field $k$ with set of even and odd generators $X=X_0 \cup X_1=\{1,x_1, \ldots, x_n\}$,
and two multiplications $\cdot$ and $\{,\}.$

We consider ideal $I$ of super-dialgebra $D_X$, generated by 
$$\big\{ \  \{a,b\}+(-1)^{|a||b|} \{b,a\}, \ a\cdot b-(-1)^{|a||b|}b\cdot a, \ 1 \cdot a-a,$$ 
$$\{a,b\}\cdot c+(-1)^{|a||b|} b\cdot \{a,c\}-D(a)\cdot b\cdot c-\{a,b\cdot c\}, $$
$$\{\{a,b\},c\} +(-1)^{|a||b|}\{b,\{a,c\}\}-\{a,\{b,c\}\} 
-D(a)\{b,c\}-(-1)^{|a||bc|}D(b)\{c,a\}-(-1)^{|c||ab|}D(c)\{a,b\} |$$ 
$$\ a,b,c \in D_{X_0} \cup  D_{X_1}   \   \big\},$$
where $D(a)=\{a,1\}.$
We say, that super-dialgebra $D_X/I$ is a free unital superalgebra of Jordan brackets.

Let $JB$ be a free unital superalgebra of Jordan brackets with set of generators $X$.
We assume, that $X$ is a linearly ordered set and $1 <x_1 < \ldots < x_n.$ 
All elements from $X$ are good words.
Let we defined all good words with length less than $n$.
The word $w$ with length $n$ in alphabet $X$ is a good word, if:

1) $w=\{u,v\}$, where $u,v$  are good words and $u>v$;

2) if $w=\{\{u_1,u_2\},v\},$ then $u_2\leq v.$

The order on good words with length less or equal $n$, we define as 
that it be consistent with order to correct for existing words, and satisfy the conditions
if word $w$ is a good and $w=\{u,v\},$ then $w>u,v.$

Define the set $M$ as
$$M=\big\{u,\{v,v\} | u, v \mbox{ are good words in alphabet  }X, v \in JB_1\big\}.$$

We define a linear order on the set $M$.
Let $u>v$, if $deg(u)>deg(v)$;
if $deg(u)=deg(v)$, 
then for $u=\{u_1,u_2\}, v=\{v_1,v_2\}$ (where $u_1>u_2, v_1>v_2$) 
we say $u>v$ if $u_1>v_1$, 
and if $u_1=v_1$ we say $u>v$ if $u_2>v_2$.
The elements of set $M$ we denote as $e_i$, where  $e_i<e_j$ if $i<j$.

Let $L_X$ be a free non-associative superalgebra over field $k$ with set of generators $X$ and multiplication $\{,\}_L$.
Consider ideal $J$ of superalgebra $L_X$, 
generated by 
$$\big\{ \  \{a,b\}+(-1)^{|a||b|} \{b,a\},$$ 
$$\{\{a,b\},c\} +(-1)^{|a||b|}\{b,\{a,c\}\}-\{a,\{b,c\}\} | \ a,b,c \in L_{X_0} \cup  L_{X_1}   \   \big\}.$$
Superalgebra $L_X/J$ is a free Lie superalgebra \cite{shtern}.
As proved in \cite{shtern}, the set $M$ is a basis of a free Lie superalgebra $L_X/J$.

Let $K$ be an ideal in $D_X$ generated by elements with type $\{ x\cdot y | x,y \neq 1\}$. 
Easy to see, that $(D_X/I)/K  \approxeq L_X/J$.

Let 
$$U=\{ e_{i_1}^{k_1} \ldots e_{i_n}^{k_n} | e_{i_1} \neq 1, e_{i_m} \in M,  i_1 < \ldots < i_n \} \cup \{1\}.$$
The set $U$ is closed under multiplication $\cdot$. 
The multiplication $\cdot$ is a supercommutative,
and we can consider elements of $M$ as 
$e_{i_1}^{k_1} \ldots e_{i_n}^{k_n},$ where $k_i \neq 0,$ and if $e_{i_j} \in JB_1$, then $k_j=1.$
We say that element $e_{i_1}^{k_1} \ldots e_{i_n}^{k_n}$ has degree $deg_e=\sum k_i$.
For element $e_i \in M$ we define $deg_L$ as length of word with generators $X$.
The degree of element $f$ is a number of generator elements in word $f$.

The main goal of this chapter is a description of basis of free superalgebra of Jordan brackets. 
We will prove that set $U$ is a basis of free superalgebra of Jordan brackets.

\begin{Lem}
The set $U$ is a basis of free superalgebra $JB$ as a vector space.
\end{Lem}

{\bf Proof.}
Similar to proof of Lemma \ref{vectorspace}.

\medskip

We define linear mapping $^*:U \rightarrow U$ as:

1) if $w \in X$, then $w^*=w$;

2) if $w \in M$, then $w^*=w$;

3) if $a=\Pi_{k=1}^n e_k^{a_k}, b=\Pi_{k=1}^n e_k^{b_k}$, then
$$(ab)^*= 
(-1)^{\sum_{i=1}^n |e_i^{a_i}| |e_{i+1}^{b_{i+1}}\cdots e_n^{b_n}| } \Pi_{k=1}^n e_k^{a_k+b_k}.$$

4) if we defined mapping $^*$ for all words with length less than $n$ and for all elements $u$ with $u<w$.
Then word $w=\{w_1,w_2\}$ 
with $w_1=\{w_{11},w_{12}\}$, 
where  $w_1,w_2,w_{11},w_{12} \in M$ and  $\{w_{11},w_{12}\}=w_1,w_{11},w_{12}>w_2,$
we can consider as element of superalgebra $(M, \{,\})$ 
and 
using (\ref{KMtojd2}) we define $w^*$ as a linear combination 
of elements from
$$W=
\big\{ \{\{w_{11},w_2\},w_{12}\}, \{w_{12},w_2 \} ,w_{11}\}, 
D(w_2) \cdot \{w_{11},w_{12}\}, D(w_{11}) \cdot \{w_2,w_{12}\}, D(w_{12})\cdot \{w_2,w_{11}\} \big\},$$
where either $u \in W$ and $u<w$ or $u \in W$ and $u=u_1\cdot u_2,$ where $u_1,u_2<w.$
It means, that for every $u \in W$ we can define $u^*$.

5) if we defined mapping $^*$ for all words with length less than $m$, 
then for word $w=\{w_1,w_2\}$ with length $m$,
where $w_1$ or $w_2$ is an element of $U \setminus M$, 
we can say that $w_2 \in U \setminus M$ and $w_2=\Pi_{k=1}^n e_k^{t_k}$,
then

$(\{w_1,w_2\})^*=$
$$\sum_{k=1}^n 
(-1)^{|e_1^{t_1}\cdots e_{k-1}^{t_{k-1}}||e_k^{t_k}|}
 t_k\left(\{w_1, e_k\}^* \frac{w_2}{e_k} \right)^* - 
\left( \sum_{k=1}^n t_k-1 \right) \left( D(w_1)^*w_2 \right) ^*,$$
where in the right part of expression elements with type $\{x,y\}^*$ was defined above.

\

%\newpage

\begin{Th}\label{basis2}
The set $U$ is a basis of free superalgebra $JB.$ 
\end{Th}

{\bf Proof.}
Similar to proof of Theorem \ref{basis}.

%We consider set $U$ as a basis of vector space and define two new multiplications  $*$ and $\{,\}^*$ by: $u*v=(uv)^*$ and $\{u,v\}^*=(\{u,v\})^*.$ It is a superalgebra $(U,*,\{,\}^*).$
%We want to prove that superalgebra $(U, *, \{,\}^*)$ is a superalgebra of Jordan brackets. After that we can say that $U$ is a basis of free superalgebra of Jordan brackets $JB.$

\medskip

Let $JB(n)$ be a space of multi-homogeneous components of free unital algebra of Jordan brackets $JB$ with $n+1$ generators $1,x_1, \ldots, x_n$.
Then 

\begin{Th}
 $dim(JB(n))=n\cdot n!.$
\end{Th}

{\bf Proof.}
Similar to proof of Theorema \ref{dimen}.

\

Using Theorem \ref{basis} and \ref{basis2} and \cite{KacKant07}, 
we can conclude

\begin{Th}
The unital superalgebra of Jordan brackets $(A,\cdot, \{,\})$ is  a free (or a simple) 
if and only if 
the
unital generalized Poisson superalgebra $(A,\cdot, \{,\}_D)$ is a free (or a simple). 
\end{Th}

%\newpage
\section{Analogues of Farkas' Theorem for PI algebras.}

The celebrated Amitsur-Levitsky Theorem  says that $M_k(R)$, 
the $k \times k$ matrices over a commutative ring $R$,
satisfies the identity $s_{2k}=0.$
The analogue of Amitsur-Levitsky Theorem for matrix superalgebras was proved in \cite{samoilov}.
Furthermore, by a well-known theorem of Amitsur, 
any associative PI algebra satisfies the condition that some $l$-th power of some $k$-th standard polynomial $s_k$ is zero: $(s_k)^l=0.$
Poisson PI algebras were studied in \cite{Fr1,Fr2,MPR,Ratseev13,Pozhidaev15}.
D. Farkas defined 
\begin{eqnarray*}\label{tojdest}
g=
\sum_{\sigma \in S_{m}} c_{\sigma} \{x_{\sigma(1)}, x_{\sigma(2)} \} \ldots \{ x_{\sigma(2i-1)} , x_{\sigma(2i)} \}
\end{eqnarray*}
as customary polynomials and proved 
that every Poisson PI algebra satisfies some customary identity \cite{Fr1}. 
Note that, in \cite{KSML} was proved the Theorem of Farkas for generic Poisson algebras. 
For generalized Poisson algebras and algebras of Jordan brackets, 
the analogue of the customary identities is 
\begin{eqnarray}\label{tojdest}
g_*=
\sum\limits_{i=0}^{[m/2]}\sum_{\sigma \in S_{m}} c_{\sigma,i} \langle x_{\sigma(1)}, x_{\sigma(2)}\rangle \ldots 
\langle x_{\sigma(2i-1)} , x_{\sigma(2i)} \rangle D(x_{\sigma(2i+1)})\ldots  D(x_{\sigma(m)}),
\end{eqnarray}
where  $$\langle x,y\rangle :=\{x,y\} -(D(x)y-xD(y)).$$

It is easy to see that algebra $(JB, \cdot, \langle , \rangle)$ satisfies identities
$$ \langle x, y\rangle= - \langle y, x\rangle,$$
$$ \langle x, y\cdot z\rangle=\langle x, y\rangle \cdot z + y \cdot \langle x,z\rangle,$$
$$ \langle\langle x,y\rangle , z\rangle + \langle\langle y,z\rangle , x\rangle +\langle\langle z,x\rangle , y\rangle =
D(x)\{y,z\}+D(y)\{z,x\}+D(z)\{x,y\}.$$
That is, $(JB, \cdot, \langle , \rangle)$ is not a Poisson algebra,
but it is a generic Poisson algebra (see \cite{KaySh12}). 
In \cite{shestakov09} it was proved that if Kantor's double of superalgebra $(JB, \cdot, \{,\})$ 
is a special Jordan superalgebra, then
$(JB,\cdot, \langle , \rangle)$ satisfies identity
$$\langle\langle a, b\rangle, c\rangle = −D(a)\langle b, c\rangle + (-1) ^{|a||b|} D(b) \langle a, c\rangle.$$

Using ideas from \cite{Fr1}, we will prove an analogue of Farkas' Theorem.
Namely, that every PI unital generalized Poisson algebra and PI unital algebra of Jordan brackets satisfies some special identity like (\ref{tojdest}).

\medskip
For homogeneous Lie word $f$ with letters $x_1,x_2, \ldots, x_n$ 
we define degree as associative monomial
$x_1^{deg_{x_1}(f)} \ldots x_n^{deg_{x_n}(f)}.$
Introduce additional useful notation
$$\{x_1, x_2,  \ldots, x_n \}: = \{ \ldots \{ x_1,x_2\} , \ldots, x_n\},$$

Following \cite{chibr06,chibr_lie}, it is easy to see  that:

\begin{Lem}\label{chi}
Let $L$ be a free Lie algebra with generators $1,y,x_1, \ldots, x_n$.
The basis of the space with components with degree $1^k y x_1 \ldots x_n$ 
is a set of elements with type
 $$ \{y, y_{\sigma(1)}, \ldots, y_{\sigma(n+k)} \},  $$
where $y_i=x_i, 0<i<n+1, y_i=1, n<i<n+k+1.$
\end{Lem}

\medskip

\begin{Lem}\label{lem41}
With notation as above,

$\{yz, w_1, \ldots, w_n \}=$
$$\sum \pm \{y, w_{i_1(1)}, \ldots, w_{i_1(k_1)} \} \{z, w_{i_2(1)}, \ldots, w_{i_2(k_2)} \} \{1, w_{i_3(1)}, \ldots, w_{i_3(k_3)} \} \ldots \{1, w_{i_l(1)}, \ldots, w_{i_l(k_l)}\},$$ 
summed over all possibilities for which each subscript appears exactly once 
in each summand and all sequences
$\{ i_1(1), \ldots, i_1(k_1) \}, \ldots, \{ i_l(1), \ldots, i_l(k_l)\}$  
are monotone increasing,
all elements $i_j(t)$ are different and $1 \leq i_j(t) \leq n$, 
also $k_1 +\ldots + k_l=n.$
(In the case of  $p=0$, the singleton $\{x\}$ is identified as $x$.)

\end{Lem}

{\bf Proof.}
Easy to see by induction.
                           
\medskip

\begin{Lem}\label{lem42} 
With notation as above,

$\{w_1, \ldots, w_r, yz, w_{r+1}, \ldots, w_{r+n}\}=$
\begin{eqnarray*}
\sum (  &\pm& \{w_0, y, w_{i_1(1)}, \ldots, w_{i_1(k_1)} \} \{z, w_{i_2(1)}, \ldots, w_{i_2(k_2)} \} \{1, w_{i_3(1)}, \ldots, w_{i_3(k_3)} \} \ldots \{1, w_{i_l(1)}, \ldots, w_{i_l(k_l)}\} \\
&\pm& \{ y, w_{i_1(1)}, \ldots, w_{i_1(k_1)} \} \{ w_0, z, w_{i_2(1)}, \ldots, w_{i_2(k_2)} \}  \{1, w_{i_3(1)}, \ldots, w_{i_3(k_3)} \} \ldots \{1, w_{i_l(1)}, \ldots, w_{i_l(k_l)}\} \\
&\pm& \{ y, w_{i_1(1)}, \ldots, w_{i_1(k_1)} \} \{ z, w_{i_2(1)}, \ldots, w_{i_2(k_2)} \} \{D(w_0), w_{i_3(1)}, \ldots, w_{i_3(k_3)} \} \ldots \{1, w_{i_l(1)}, \ldots, w_{i_l(k_l)}\}),
\end{eqnarray*}
where 
$w_0= \{w_1, \ldots, w_r\}.$
\end{Lem}
 
{\bf Proof.}
If $r=0$ this lemma follows immediately from the previous one. 
So we assume that $r \geq 1$ and 
for $\beta=\{yz, w_0, w_{r+1}, \ldots, w_{r+n}\}$
we apply the previous lemma formally to this new expression.
Easy calculations give the proof of the lemma.

\medskip

\begin{Lem}\label{lem43}
Assume $f$  is a multi-linear Lie polynomial in the letters $x, w_1, \ldots, w_m$, 
 and $f$ is an associative algebraic derivation in $x$ when considered as a Poisson polynomial.
Then $deg(f) = 2.$
\end{Lem}

{\bf Proof.}
Let $f$ have a degree $1^kx w_1\ldots w_m.$
Set $w_i=1, m<i<m+k+1$ and $x=w_{m+k+1}.$
By Lemma \ref{chi}, we may write
$$h=\sum c_{\sigma} \{ w_1, w_{\sigma(2)}, \ldots,  w_{\sigma(m+k+1)} \},$$
where $\sigma$ runs over permutations of $2, \ldots, m+k+1$ and $c_{\sigma}$ is a scalar.
We compute the derivation difference of $h$ with respect to $w_{m+k+1}\rightarrow yz$.
%Zoom in on the derivation difference of the 

One of the many summands will be the product of $c_{\sigma}$ and
$$\{w_1, w_{\sigma(2)}, \ldots , w_{\sigma(t-1)}, y\}\{z, w_{\sigma(t+1)}, \ldots, w_{\sigma(m+k+1)}\}.$$
And one of the many others summands will be product of $c_{\sigma}$ and $yzD\{w_1, w_{\sigma(2)}, \ldots, w_{\sigma(m+k)}\}.$
But this expressions determines $\sigma$;
it does not appear in the derivation difference for any other permutation.
If we specify an ordering of the variables 
$$y<z<w_1< \ldots < w_n,$$
then these special expressions are linearly independent for $deg(f)>2$ by the application following Lemma \ref{chi}.
%Since the derivation difference of $h$ is zero, the coeff
The Lemma is proved.

\medskip

For polynomial $f(x, \ldots)$ we define  
$$\Delta_{yz}^{x}(f)  := f(yz, \ldots)-yf(z, \ldots)-z f(y, \ldots).$$

It is easy to see, that if $\Delta_{yz}^{x}(f)\neq 0$
and algebra $A$ satisfies identity $f$, then $A$ also satisfies identity $\Delta_{yz}^{x}(f)$.

For Poisson polynomial $f$ we define $x$-$height(f)$ as maximal length of $\{,\}$-subword with $x$.

\medskip

\begin{Lem}\label{lem44}
If multi-linear Poisson polynomial $f$ is an associative derivation in any letter $x$, 
then $x$-$height(f)<3.$
\end{Lem}

{\bf Proof.}
%Предположим, что существует слово $g$ степени $xm+p$, где $deg(m)>1$, причем $deg(m)$ максимальна для $f$.
Let $f=\sum_{\gamma} f_{\gamma}\gamma,$ where $f_{\gamma}$ be a Lie polynomials with degree $1^kxt$ for any $t$. 
We consider $f_{\gamma}=\sum_k f_{\gamma}^k,$ where the degree of polynomial $f_{\gamma}^k$ is $1^k xt.$
On polynomials $f_{\gamma}^k$ we choose a polynomials with maximum $deg(t)$ 
and after that polynomial with the maximum $k$.
Let it is $f^*$ with degree $1^kxx_1\ldots x_n.$
Continuing the argument similar to that of Lemma \ref{lem43},
we obtain  that if $deg(f^*)>2,$ then $\Delta_{yz}^{x}(f) \neq 0.$
%По лемме \ref{chi}, мы можем представить  $$f^*= \sum_{\sigma} c_{\sigma} \{x, y_{\sigma(1)}, \ldots, y_{\sigma(n+k)} \},$$ где $y_i=x_i, 0<i<n+1, y_i=1, n<i<n+k+1.$
%Откуда, видим, что  $f(x \rightarrow yz)- yf(x \rightarrow z)-f(x \rightarrow y)z$ имеет слагаемые вида $\{y, w_*\}\{z,w_*\}\gamma,$ которые линейно независимы c другими слагаемыми в силу Lemma \ref{chi}.
%Using Lemma \ref{lem43}, $deg(f^*)<3.$
Following, $deg(f_{\gamma})<3$ and $x$-$height(f)<3.$
The Lemma is proved.

\medskip

The main result of this chapter is following theorem.

\medskip

\begin{Th}\label{ThPI}
If unital generalized Poisson algebra $A$ satisfies polynomial identity $g_0$, 
then $A$ satisfies polynomial identity $g_*$ of type (\ref{tojdest}).
\end{Th}

{\bf Proof.} Using ideas from \cite{Fr1}, we consider three steps.

{\bf Step 1.}
Let $x$-$height(g_0)>2.$
Then algebra $A$ 
satisfies identity $g_1=\Delta_{yz}^{x}(g_0).$ By Lemma \ref{lem44}, $g_1\neq 0$.
Note that $x$-$height(g_0)>y$-$height(g_1),z$-$height(g_1).$
If polynomial $g_1$ has letter $w$ with $w$-$height(g_1)>2$, then we will continue the process.
That is, 
%Таким образом, в силу монотонной убываемости значений $w$-$height$ по каждой букве $w$,
in the end we can obtain that algebra $A$ satisfies some identity $f$, 
where arbitrary letter $v$ has $v$-$height(f)<3.$

{\bf Step 2.}
Let algebra $A$ satisfy identity $f$ and $x$-$height(f)<3$ for arbitrary letter $x$ from $f$.
Assume that $f$ is not an associative derivation for any letter $x$. 
Then $$f=xT+D(x)T_0+\sum \{x,x_i\} T_i.$$
We consider $$\Delta^x_{yz}(f)= -yzT +yz \sum D(x_i)T_i \neq 0.$$ 
It is easy to see that 
unital algebra $A$ will satisfy identity $f^*=-T + \sum D(x_i)T_i,$ 
with less letters than identity $f$. 
That is, 
if identity $f^*$ 
is not a derivation for any letter $t$, then we can continue this process. 
The process is finite,
and algebra $A$ satisfies some identity $g$ with $x$-$height(g)<3$,
and $g$ is an associative derivation on every letter. 

{\bf Step 3.}
For polynomial $g$ with letters $x_1, \ldots, x_n,$ 
every associative subword with type $\{x,y\}$ 
we can consider as $\langle x,y\rangle +(D(x)y-xD(y)).$
It follows that polynomial $g$ we can consider as $g_*+D(x_1)T_1+x_1R_1$, where $g_*$ has type (\ref{tojdest}).
It is easy to see that $g_*$ is an associative derivation on every letter, 
then algebra $A$ will satisfy identity $R_1.$
It follows that $A$ satisfies identity $g_*+D(x_1)T_1$. 
We will consider $D(x_1)T_1$ as $D(x_2)T_2 +x_2R_2$.
That is, algebra $A$ satisfies an identity $g_*+D(x_1)D(x_2)T_{1,2}$. 
Continuing this process, 
we find that algebra $A$ satisfies identity $g_*$ of the form (\ref{tojdest}).
The Theorem is proved.

\medskip

As follows, we can obtain

\medskip

\begin{cor}\label{cor}
Every PI generalized Poisson algebra satisfies an identity of the following type
\begin{eqnarray}\label{pi}
 f_{*} =
\sum\limits_{i=0}^{[m/2]}\sum_{\sigma \in S_{m}} c_{\sigma,i}  
 \prod\limits_{k=1}^{i} [x_{\sigma(2k-1)}; x_{\sigma(2k)}; z_{2k-1}; z_{2k}]  \cdot 
 \prod\limits_{k=1}^{m-2i} \{x_{\sigma(2i+k)}; z_{2i+2k-1}; z_{2i+2k} \}  \cdot 
 \prod\limits_{k=1}^{2i} z_{2m-2i+k},
\end{eqnarray}
where 

$[u_1; u_2; w_1;w_2]=$
$$\{u_1,u_2\}w_1w_2+\{u_1,w_1w_2\}u_2 + u_1\{w_1w_2, u_2\}+
\sum _{\sigma_1, \sigma_2 \in S_2} \{u_{\sigma_1(1)}, w_{\sigma_2(1)}\}u_{\sigma_1(2)}w_{\sigma_2(2)},$$

$$\{t_1;t_2;t_3 \} = \{t_2t_3,t_1\} - \{t_2,t_1\}t_3-\{t_3,t_1\}t_2.$$
\end{cor}

{\bf Proof.}
By theorem \ref{ThPI}, every PI generalized Poisson algebra satisfies an  identity $g_*$ of type (\ref{tojdest}).
Multiplying $g_*$ and $\Pi_{k=1}^{2m} z_k$,
note that
$$w_1 w_2 \langle u_1,u_2 \rangle = [u_1; u_2; w_1;w_2]  \mbox{ and }t_2t_3 D(t_1)=\{t_1;t_2;t_3 \}.$$
It is easy to see the corollary is proved.

\medskip

Also, we can prove an analogue of Farkas' Theorem for PI algebras of Jordan brackets:

\begin{Th}\label{ThPI2}
If unital algebra of Jordan brackets $A$ satisfies polynomial identity $g_0$, 
then $A$ satisfies polynomial identities of types (\ref{tojdest}) and (\ref{pi}).
\end{Th}

{\bf Proof.}
We consider generalized Poisson algebra $A_D$ with vector space $A$ and multiplications $\cdot$ and $\{,\}_D$.
It is easy to see, that algebra $A_D$ is a PI algebra.
Using Theorem \ref{ThPI}, we can obtain that algebra $A_D$ satisfies a polynomial identity $g$ of type (\ref{tojdest}).
Noted that, $2D_D=D$ and
$$\{x,y\}_D+xD_D(y)-D_D(x)y=\{x,y\}+xD(y)-D(x)y,$$
we conclude that algebra $A$ satisfies a polynomial identity of type (\ref{tojdest}). 
And using proof of Corollary \ref{cor}, we find that algebra $A$ satisfies polynomial identity of type (\ref{pi}).
The Theorem is proved.

\section{Kantor's Double for generic Poisson superalgebras.}

Let us give the definition of a generic Poisson superalgebra \cite{KSML}. 
We consider a vector space $GP$ with associative supercommutative multiplication $x\cdot y$
and superanticommutative multiplication $\{x,y\}.$ If superalgebra $(GP,\cdot,\{,\})$
satisfies the Leibniz identity, i.e., \begin{eqnarray}
\{x, y \cdot z\}=\{x,y\}\cdot z+(-1)^{|x||y|}y\cdot\{x,z\},\label{gpident}\end{eqnarray}
then $GP$ is a generic Poisson superalgebra (GP superalgebra). 
It is easy to see that if $(GP,\{,\})$ is a Lie superalgebra, then $(GP,\cdot,\{,\})$ is a Poisson superalgebra.
It is well-known that Kantor's Double of a Poisson superalgebra is a Jordan superalgebra \cite{Kantor90}.
There naturally arises a \\

{\bf Question.}
When is Kantor's double of a generic Poisson superalgebra  a Jordan superalgebra?

\

We gave the answer on this question in following\\

\begin{Th}
Let $A$ be a GP superalgebra over field $F$ with characteristic $\neq 2,3$, 
then $K(A)$ is a Jordan superalgebra if and only if $A$ satisfies special identity
\begin{eqnarray}\label{jordanGP}
(\{\{a,b\},c\}-(-1)^{|b||c|}\{\{a,c\},b\} - \{a,\{b,c\}\})d=0
\end{eqnarray}
for homogeneous elements $a,b,c,d \in GP_0 \cup GP_1.$
\end{Th}

\

{\bf Proof.}
In \cite{Kay10} it was proved that
if $(B, \cdot)$ is an associative supercommutative superalgebra and $(B, \{,\})$ is a superanticommutative superalgebra,
then $K(B)$ is a Jordan superalgebra
if and only if for homogeneous elements superalgebra $B$ satisfies following identities

$$(-1)^{(k+j)i}(\{h_kk_l,g_j\}f_i-h_kk_l\{g_j,f_i\})=$$
\begin{eqnarray}\label{jorskob2}(-1)^{(l+j)k}(\{k_lf_i,g_j\}h_k-k_lf_i\{g_j,h_k\}),\end{eqnarray}

$$(-1)^{(i+j)l}(\{f_ih_kg_j,k_l\}-f_ih_k\{g_j,k_l\})=$$
$$(-1)^{(k+j)i}(\{h_kk_l,g_j\}f_i-\{h_kk_l,g_jf_i\})+$$
\begin{eqnarray}\label{jorskob3}
(-1)^{(l+j)k}(\{k_lf_i,g_j\}h_k-\{k_lf_i,g_jh_k\}),\end{eqnarray}

$$(-1)^{(i+j)l}\{\{f_{i},h_k\}g_j,k_l\}+(-1)^{(k+j)i}\{\{h_k,k_l\}g_j,f_i\}+$$
$$(-1)^{(l+j)k}\{\{k_l,f_i\}g_j,h_k\}=(-1)^{(i+j)l}\{f_i,h_k\}\{g_j,k_l\}+$$
\begin{eqnarray}\label{jorskob1}(-1)^{(k+j)i}\{h_k,k_l\}\{g_j,f_i\}+(-1)^{(l+j)k}\{k_l,f_i\}\{g_j,h_k\},\end{eqnarray}
where $k_i,h_i, g_i,f_i \in B_i.$

It is easy to see that if an algebra satisfies identity (\ref{gpident}),
then it satisfies identities (\ref{jorskob2}, \ref{jorskob3}).
Using identity (\ref{gpident}) for (\ref{jorskob1}),
we get (\ref{jordanGP}). The Theorem is proved.

\begin{cor}
If $A$ is a unital GP (non-Poisson) superalgebra, then $K(A)$ is not a Jordan superalgebra.
\end{cor}

\begin{cor}
If $A$ is a GP (non-Poisson) superalgebra and $(A,\cdot)$ is an algebra with trivial annihilator (for example, prime algebra), 
then $K(A)$ is not Jordan superalgebra.
\end{cor}

\begin{exa}
Let $A$ be a non-Lie anticommutative algebra with multiplication $\{ ,\}$ (for example, a Malcev  or  a binary-Lie algebra).
We define new multiplication $x\cdot y =0$.
Then, algebra $(A, \cdot, \{ ,\})$ is a GP (non-Poisson) algebra, 
and satisfies identity (\ref{jordanGP}),
and $K(A)$ is a Jordan superalgebra.
\end{exa}

{\bf Acknowledgments:}
I am grateful to Prof. Dr. Ivan Shestakov (IME-USP, Brazil) for the idea of this work and some interesting comments.

\

%\newpage

\end{document}